\documentclass[11pt]{article}          
\usepackage{amstex}                      
\usepackage[psamsfonts]{amssymb}    

\begin{document}

\def\claim#1{\noindent {\bf{#1}}\ \ }
\def\dim{{\rm dim}}         
\def\K{{\cal K}}         
\def\G{{\cal G}}         
\def\b{\beta}         
\def\R{{\mathbb R}}     
\def\Q{{\mathbb Q}}     
\def\C{{\mathbb C}}     
\def\Z{{\mathbb Z}}     
\def\N{{\mathbb N}}     
\newcommand{\BoX}{{}\raisebox{-.18ex}{\makebox[1em][r]{$\Box$}}}
\newcommand{\blok}{\hspace*{0cm}\hfill\BoX\par}


\newtheorem{theorem}{Theorem}[section]
\newtheorem{lemma}[theorem]{Lemma}
\newtheorem{proposition}[theorem]{Proposition}
\newtheorem{definition}[theorem]{Definition}
\newtheorem{corollary}[theorem]{Corollary}
\newtheorem{fig}{Figure}

\centerline{{\LARGE \bf Vassiliev Invariants and Gleam Polynomials}}
\medskip
\bigskip

\bigskip
\centerline{Urs Burri}
\bigskip

\centerline{Preprint}

\begin{abstract}
\noindent
It was shown by Goussarov~\cite{gu1} that Vassiliev invariants 
are polynomials in the gleams for a fixed Turaev shadow. In this 
paper we show that Vassiliev invariants are almost characterized 
by this fact. We also prove that the space of knot invariants 
which are polynomials in the gleams is bigger than the Vassiliev 
subspace.
\end{abstract}

\section*{Introduction}                    \label{introduction}

In Turaev's shadow topology \cite{tu} a knot in $S^3$ is 
presented by its Turaev shadow, which is its image under 
the Hopf projection, and its gleams, which are integer 
numbers associated to the regions which enable us to 
reconstruct the knot up to isotopy. It is a very interesting 
question how knot invariants depend on the gleams for a fixed 
Turaev shadow. Studying this question or searching for explicit
 formulas it is easier to work with the parametrization $\K$
 introduced in \cite{bu}, than with the gleams directly. This
 map $\K$ is defined on a lattice $\Z^e\times \Z$, where $e$
 is the number of double points of the fixed Turaev shadow. 
One of the first results was the following: Let us define 
Jones' Vassiliev invariants $\{u_n\}_{n\geq 2}$ via the 
Jones polynomial $J_t(K)$, by 
$J_{e^x}(K)=\sum_{n=0}^{\infty}u_n(K)x^n$.

\begin{theorem}{\em \cite{bu}}               \label{bu}
The function $u_n\circ \K :\Z^e \times \Z\longrightarrow \Q$ 
is a polynomial of degree $2n$.
\end{theorem}
Before we state the main results, let us recall the following.
 A knot invariant $v$ with values in some abelian group $\G$ 
is a polynomial of degree $m$ in the gleams, if and only if
 $v\circ \K$ is a polynomial of degree $m$.

\section{The Main Results}                    \label{results}

\noindent

\noindent
In the following we will state all results for knots only, 
although most of them easily extend to links. Let $v_n$ be
 a Vassiliev invariant of order $n$ with values in $\G$.

\begin{theorem}   {\em (Goussarov \cite{gu1})}       \label{gu1}
The function $v_n\circ \K :\Z^e \times \Z\longrightarrow \G$
 is a polynomial of degree $\leq 2n$.
\end{theorem}

\noindent
\claim{Remark} Goussarov proved that derivatives of order 
$2n+1$ vanish. However, he cannot give explicit formulas
 for $v_n\circ \K $.

\begin{theorem}                                 \label{gu2}
Let $v$ be a knot invariant with values in $\G$. Assume 
that for every fixed Turaev shadow $v$ is a polynomial 
in the gleams of degree $\leq m$. Then $v$ is a Vassiliev 
invariant of order $\leq m$. 
\end{theorem}

\noindent
\claim{Remark} The main point (due to Goussarov \cite{gu2}) 
is that one should not try to prove that a knot invariant
 which is a polynomial in the gleams of degree $\leq 2m$ 
is a Vassiliev invariant of order $\leq m$.
\\

\noindent
{\it Proof.} We have to show that $v$ vanishes after the usual
 extension to singular knots on every $(m+1)$-singular knot
 $K^{(m+1)}$. Suppose that $K^{(m+1)}$ is given by a standard
 diagram on $\R^2$. We would like to apply Thereom 3.2 in 
\cite{bu} where alternating sums of values of a knot invariant
 $v$ are related with partial derivatives (finite differences)
 of $v\circ \K$. First, we fix the Turaev shadow $s$ on 
$\R^2\cup \{\infty\}$ given by the projection of $K^{(m+1)}$.
 We enumerate the double points of $s$. Let $P_1,...,P_{m+1}$ 
be the double points of $s$ coming from the $m+1$ singular 
points of $K^{(m+1)}$ and $P_{m+2},...,P_e$ be the remaining
 double points of $s$ coming from the negative and positive
 crossings of $K^{(m+1)}$. Then we know that
$$v\left(K^{(m+1)}\right)=\Delta_{x_1}\Delta_{x_2}...
\Delta_{x_{m+1}}(v\circ \K)(0,0,...,0,
\b_{m+2},...,\b_{e},0)\ ,$$
where $\b_k$ is $0$ or $1$ if the crossing in $P_k$ of 
$K^{(m+1)}$ is negative or positive, respectively. However,
 by assumption, $v\circ\K$ is a polynomial of degree 
$\leq m$, so the $m+1$ partial derivatives of 
$v\circ \K$ vanish and we are done. \blok

\begin{theorem}                                    \label{gu3}
There exist $\C$-valued knot invariants, which are 
polynomials in the gleams for every fixed Turaev 
shadow, but are not Vassiliev.
\end{theorem}

\noindent
\claim{Remark} The point here is that the degree of these 
polynomials is not uniformly bounded, when we change the
 Turaev shadows. The idea of the following construction
 is due to Goussarov \cite{gu2}.
\\

\noindent
{\it Proof.} We use the fact, that the dimensions 
$\dim_{\C}(V_n)$ of the vector spaces of Vassiliev 
invariants of order $\leq n$ grow faster than 
polynomially, when $n$ tends to $+\infty$, see 
\cite{ko}. Kontsevich states that
 $$\dim_{\C} (V_n)>e^{c\sqrt{n}} \ , 
\ \quad n \to +\infty \ ,$$
 for any positive constant $c<\pi \sqrt{2/3}$. Let us denote
 by $m_k$ the number of Turaev shadows with $\leq k$ regions,
 and by $W_{2n,k}$ the vector space of polynomials 
(not necessarily knot invariants) of degree 
$\leq 2n$ in $\leq k$ variables. We define by 
$\pi_{n,k}$ the linear map which associates to 
any Vassiliev invariant of order $\leq n$, 
the $m_k$ polynomials of order $\leq 2n$ in $\leq k$ 
variables, see Theorem~\ref{gu1}.
$$\pi_{n,k}:V_n\longrightarrow (W_{2n,k})^{\oplus m_k}\ .$$ 
The dimension of $V_n/{V_{n-1}}$ grows faster in $n$ 
than $\dim (W_{2n,k})^{\oplus m_k}$ for any fixed $k$;
 the first growth is subexponential, the second one
 polynomial. Therefore there exists for any $k$ a 
Vassiliev invariant $v_{n_k} \in V_{n_k} \setminus
 V_{n_k-1}$ with $\pi_{n_k,k}(v_{n_k})=0\in 
(W_{2n_k,k})^{\oplus m_k}$. This means that 
$v_{n_k}$ is a Vassiliev invariant of order 
$n_k$ (not of any lower order), which 
vanishes on all knots which admit a Hopf
 projection with $\leq k$ regions. Now define 
$$v:=\sum_{k=0}^{\infty}v_{n_k}.$$
Note that $v$ is well defined because for every 
knot the sum is finite. \blok

\medskip

\noindent
\claim{Remark} This result shows that the space 
of knot invariants which are polynomials in the 
gleams is bigger than the Vassiliev subspace. It 
is known that many ``classical'' knot invariants
 are not Vassiliev invariants but unfortunately 
they are not polynomials in the gleams either, as 
we will see in the next theorem. Let us mention 
one more thing first. In \cite{tr} it was shown
 that the uniform limit of Vassiliev invariants 
is again Vassiliev. Using the fact that these 
Vassiliev invariants give rise to polynomials
 defined on the lattice $\Z^e\times \Z$ it is
 easy to check that such a sequence
 $(v^{(n)})_{n\in \N}$ has to stabilize, 
if it is normalized, let us say by requiring
 $v^{(n)}({\rm unknot})=0$. The same is true 
for a uniform convergent sequence of knot 
invariants, which are polynomials in the gleams.

\begin{theorem}                             \label{not}
The following knot invariants are not polynomials in 
the gleams: the unknotting number, genus, signature,
 bridge number, braid index, span of the Jones polynomial.
\end{theorem}

\noindent
{\it Proof.} It is enough to study the behaviour of 
these invariants on the series of knots denoted by
 $\K(x_f)$, see \cite{bu}. These are the only 
knots which admit a Hopf projection without any
 double points. Recall that for $x_f \geq 2$,
 $\K(x_f)={\rm torus}\ {\rm knot}(x_f,x_f+1)$
 and $\K(x_f)$ is trivial for $x_f \in \{-2,-1,0,1\}$.
 All the knot invariants mentioned above have the same
 behaviour: They are not constant on $\K(x_f)$ but 
they grow only linearly or at most quadratically 
in $x_f$, and therefore they are not polynomials in $x_f$.
\blok


\def\papername{\it}

\bigskip
\bigskip
\noindent
Urs Burri\par
\noindent
Mathematisches Institut Basel/Bern\par
\noindent
Switzerland\par
\noindent
burri@@math.unibas.ch

\end{document}